\documentclass[12pt,a4paper]{amsart}
\usepackage{amssymb}
\usepackage{mathrsfs}
\headsep=1cm \numberwithin{equation}{section}
\newtheorem{thm}{Theorem}[section]
\newtheorem{cor}[thm]{Corollary}
\newtheorem{lem}[thm]{Lemma}

\newtheorem{prop}[thm]{Proposition}
\theoremstyle{definition}
\newtheorem{defn}[thm]{Definition}

\newtheorem{rem}[thm]{Remark}
\numberwithin{equation}{section}

\newcommand\Supp{\operatorname{Supp}}
\newcommand\Ass{\operatorname{Ass}}
\newcommand\mAss{\operatorname{mAss}}

\newcommand\Ann{\operatorname{Ann}}
\newcommand\Spec{\operatorname{Spec}}
\newcommand\Rad{\operatorname{Rad}}

\newcommand\Ext{\operatorname{Ext}}

\newcommand\Att{\operatorname{Att}}

\begin{document}\title[Annihilators and attached primes of local cohomology]
{On the annihilators and attached primes of top local cohomology modules}
\author [Ali Atazadeh, Monireh Sedghi and Reza Naghipour]{Ali Atazadeh, Monireh Sedghi$^*$ and Reza Naghipour\\\\\\\,
\vspace*{0.5cm}Dedicated to Professor Hossein Zakeri}
\address{Department of Mathematics, Azarbaijan Shahid Madani University, Tabriz, Iran.}
\email{aalzp2002@yahoo.com}

\address{Department of Mathematics, Azarbaijan Shahid Madani University, Tabriz, Iran.}
\email {m\_sedghi@tabrizu.ac.ir}
\email {sedghi@azaruniv.ac.ir}

\address{Department of Mathematics, University of Tabriz, Tabriz, Iran;
and School of Mathematics, Institute for Research in Fundamental
Sciences (IPM), P.O. Box 19395-5746, Tehran, Iran.}
\email{naghipour@ipm.ir} \email {naghipour@tabrizu.ac.ir}
\thanks{ 
The research of the third author was supported in part by a grant from IPM.\\
$^*$Corresponding author: e-mail: {\it m\_sedghi@tabrizu.ac.ir and sedghi@azaruniv.ac.ir} (Monireh Sedghi)}%

\begin{abstract}
Let $\frak a$ be an ideal of a commutative Noetherian ring $R$ and $M$ a finitely generated $R$-module.  It is shown that $\Ann_R(H_{\frak a}^{\dim M}(M))= \Ann_R(M/T_R(\frak a, M))$, where  $T_R(\frak a, M)$ is the largest submodule of $M$ such that ${\rm cd}(\frak a, T_R(\frak a, M))< {\rm cd}(\frak a, M)$. Several applications of this result are given. Among other things, it is shown that  there exists an ideal $\frak b$ of $R$ such that $\Ann_R(H_{\frak a}^{\dim M}(M))=\Ann_R(M/H_{\frak b}^{0}(M))$. Using this, we  show that if $ H_{\frak a}^{\dim R}(R)=0$, then
$\Att_RH^{\dim R-1}_{\frak a}(R)= \{\frak p\in \Spec R|\,{\rm cd}(\frak a, R/\frak p)=\dim R-1\}.$
These generalize the main results of \cite[Theorem 2.6]{BAG}, \cite[Theorem 2.3]{He} and \cite[Theorem 2.4]{Lyn}.
\end{abstract}
\maketitle
\section{Introduction}
Let $R$ be an arbitrary commutative Noetherian ring
(with identity), $\frak a$ an ideal of $R$ and let $M$ be a finitely generated $R$-module. An important problem concerning local cohomology is determining the annihilators of the $i^{\rm th}$ local cohomology module $H_{\frak a}^i(M)$. This problem has been studied by several authors; see for example \cite{HK}, \cite{Lyn}, \cite{Ly1}, \cite{NA},  \cite{NC},  \cite{Sc1},   and  has led to some interesting results. More recently,
in \cite{BAG} Bahmanpour et al., proved an interesting result about the annihilator $\Ann_R(H_{\frak m}^d(M))$,  in the case $(R, \frak m)$ is a complete local ring of dimension $d$.

The purpose of the present paper  is to establish some new results concerning of the annihilators of local cohomology modules $H_{\frak a}^i(M)$ $(i\in \mathbb{N}_0)$, where $\frak a$ is an ideal in a Noetherian ring $R$ and $M$ a finitely generated module over $R$.

As a main result in the second section, we determine the  annihilators of the top local cohomology module
$H_{\frak a}^{\dim M}(M)$.   More precisely, we shall prove the following theorem which is a generalization of the main result of \cite[Theorem 2.6]{BAG} for an arbitrary ideal $\frak a$ of an  arbitrary  Noetherian ring $R$.

\begin{thm}
Let $R$ be a Noetherian ring and let $M$ be a finitely
generated $R$-module. Then for any ideal $\frak a$ of $R$, $\Ann_R(H_{\frak a}^{\dim M}(M))= \Ann_R(M/T_R(\frak a, M))$,
where $T_R(\frak a, M)$ denotes the largest submodule of $M$ such that ${\rm cd}(\frak a, T_R(\frak a, M))< {\rm cd}(\frak a, M)$.
\end{thm}

The result  in Theorem 1.1 is proved in Theorem 2.3.
Several corollaries of this result are given.  A typical result in this direction is the following, which is a generalization of the main results of \cite[Theorem 2.6]{BAG} and \cite[Theorem 2.4]{Lyn}  for an  ideal $\frak a$ in an arbitrary Noetherian ring $R$.
\begin{cor}
Let $R$ be a Noetherian  ring and $\frak a$ an ideal of $R$. Let $M$ be a non-zero finitely generated $R$-module of finite  dimension $c$ such that
${\rm cd}(\frak a, M)=c$. Then $$\Ann_R(H_{\frak a}^{c}(M))=\Ann_R(M/H_{\frak b}^{0}(M))=\Ann_R(M/\cap_{{\rm cd}(\frak a, R/{\frak p_j})=c}N_j).$$ Here $0=\cap_{j=1}^nN_j$ denotes a reduced primary decomposition of the zero submodule $0$ in $M$ and $N_j$ is a $\frak p_j$-primary submodule of $M$, for all $j=1,\dots ,n$ and $\frak b:=\Pi_{{\rm cd}(\frak a, R/{\frak p_j})\neq c}\frak p_j.$
\end{cor}

One of the basic problems concerning local cohomology is finding the set of attached primes of  $H_{\frak a}^{i}(M)$.
In  Section 3, we obtain some results about the attached primes of  local cohomology modules. In this section among other things, we derive the following consequence of Theorem 1.1 and Corollary 1.2, which provides an upper bound for the attached primes of
$\Att_RH^{{\rm cd}(\frak a, M)}_{\frak a}(M)$.  This will generalize the main results of \cite{Di} and \cite{DY}.

\begin{thm}
Let $R$ be a Noetherian ring and $\frak a$ an ideal of $R$. Let $M$ be a finitely generated $R$-module such that $c:={\rm cd}(\frak a, M)$ is finite. Then
$$ \Att_RH^c_{\frak a}(M)\subseteq \{\frak p\in \Supp M|\,{\rm cd}(\frak a, R/\frak p)= c\}.$$
Moreover, if $c=\dim M$, then  $$ \Att_RH^c_{\frak a}(M)= \{\frak p\in \mAss_R M|\,{\rm cd}(\frak a, R/\frak p)= c\}.$$
\end{thm}

For an $R$-module $A$, a prime ideal $\frak p$ of $R$  is said to be {\it attached prime to} $A$ if $\frak p=\Ann_R(A/B)$ for some submodule $B$
of $A$. We denote the set of attached primes of $A$ by $\Att_RA$.  This definition agrees with the usual definition of  attached  prime if $A$ has a secondary representation (cf. \cite[Theorem 2.5]{Ma}).

Another main result in Section 3 is to give a complete characterization of the attached primes of the local cohomology module $H_{\frak a}^{\dim R-1}(R)$. More precisely, we shall show the following result,  which is an extension, as well as, a correction of the main theorem of \cite{He}.\\
\begin{thm}
Let $(R, \frak m)$ be a local (Noetherian) ring of dimension $d$. Let $\frak a$ be an ideal of $R$ such that $\dim R/\frak a=1$ and
$H^d_{\frak a}(R)=0$.  Then ${\rm Assh}_RR\subseteq \Att_RH^{d-1}_{\frak a}(R).$ Moreover, if $R$ is complete, then
$$\Att_RH^{d-1}_{\frak a}(R)=\{\frak p\in \Spec R|\, \dim R/\frak p=d-1\,\text{and}\, \Rad(\frak a+\frak p)=\frak m\}\cup {\rm Assh}_RR.$$
\end{thm}

One of our tools for proving Theorem 1.4 is the following:\\

\begin{prop}
 Let $R$ be a Noetherian ring of finite dimension $d$ and $\frak a$ an ideal of $R$ such that $H^d_{\frak a}(R)=0$. Then
$$\Att_RH^{d-1}_{\frak a}(R)= \{\frak p\in \Spec R|\,{\rm cd}(\frak a, R/\frak p)=d-1\}.$$
\end{prop}

Throughout this paper, $R$ will always be a commutative Noetherian
ring with non-zero identity and $\frak a$ will be an ideal of $R$. For any $R$-module $M$, the
$i^{\rm th}$ local cohomology module of $M$ with support in $V(\frak a)$
is defined by
$$H^i_{\frak a}(M) := \underset{n\geq1} {\varinjlim}\,\, \Ext^i_R(R/\frak a^n, M).$$

 The cohomological dimension of $M$ with respect to $\frak a$ is defined as $${\rm cd}(\frak a, M):= \sup\{i\in \mathbb{Z}|\, H^i_{\frak a}(M)\neq 0\}.$$

For each $R$-module $L$, we denote by
 ${\rm Assh}_RL$ (resp. $\mAss_RL$) the set $\{\frak p\in \Ass
_RL:\, \dim R/\frak p= \dim L\}$ (resp. the set of minimal primes of
$\Ass_RL$).  Also,  we shall use $\Att_R L$ to denote the set of attached prime
ideals of $L$.
For any ideal $\frak a$ of $R$, we denote $\{\frak p \in {\rm Spec}\,R:\, \frak p\supseteq \frak a \}$ by
$V(\frak a)$. Finally, for any ideal $\frak{b}$ of $R$, {\it the radical} of $\frak{b}$, denoted by $\Rad(\frak{b})$, is defined to
be the set $\{x\in R \,: \, x^n \in \frak{b}$ for some $n \in \mathbb{N}\}$.  For any unexplained notation and terminology we refer
the reader to \cite{BS} and \cite{Mat}.

\section{Annihilators of top local cohomology modules}

Let us, firstly, recall the important concept of the cohomological dimension of an $R$-module $L$ with respect to an ideal $\frak a$ of a commutative Noetherian ring $R$. Denoted by ${\rm cd}(\frak a,L)$, is the largest integer $i$ such that $H^i_{\frak a}(L)\neq 0$. The main result of this section is Theorem 2.3. The following lemma plays a key role in the proof of that theorem.\\

\begin{lem}
Let $R$ be a Noetherian ring and $\frak a$ an ideal of $R$. Let $M$ and $N$ be two finitely generated
$R$-modules such that $\Supp N\subseteq \Supp M$. Then, ${\rm cd}(\frak a, N)\leq {\rm cd}(\frak a, M).$
\end{lem}
\proof See \cite[Theorem 2.2]{DNT}. \qed\\

\begin{defn}
Let $R$ be a Noetherian ring and $\frak a$ an ideal of $R$. Let $M$ be a non-zero
finitely generated $R$-module. We denote by ${\rm T}_R(\frak a, M)$ the largest submodule of $M$ such that ${\rm cd}(\frak a, {\rm T}_R(\frak a, M))<{\rm cd}(\frak a, M)$. It is easy to check that ${\rm T}_R(\frak a, M)=\cup \{N|\,N\leq M$ and ${\rm cd} (\frak a, N)<{\rm cd} (\frak a, M)\}$. In particular for a local ring $(R, \frak m)$, we denote ${\rm T}_R(\frak m, M)$ by $T_R(M)$. Thus $T_R(M)=\cup \{N|\,N\leq M$ and $\dim N<\dim M\}$, see \cite[Definition 2.5]{BAG}.\\
\end{defn}
Now, we are prepared to present the main result of this section which is a generalization of the main result of \cite[Theorem 2.6]{BAG}.\\

\begin{thm}
Let $R$ be a Noetherian ring and $\frak a$  an ideal of $R$. Let $M$ be a finitely generated
$R$-module with finite dimension $c$ such that ${\rm cd}(\frak a, M)=c$. Then $$\Ann_R(H_{\frak a}^c(M))= \Ann_R(M/T_R(\frak a, M)).$$
\end{thm}
\proof It easily follows from the canonical exact sequence $$0 \longrightarrow T_R(\frak a, M) \longrightarrow M \longrightarrow M/T_R(\frak a, M) \longrightarrow 0,$$ and Lemma {\rm2.1} that ${\rm cd}(\frak a, M/T_R(\frak a, M))=c$ and $$H_{\frak a}^c(M)\cong H_{\frak a}^c(M/T_R(\frak a, M)).$$ Thus we have $$H_{\frak a}^c(M)\cong H_{\frak a}^{{\rm cd}(\frak a, M/T_R(\frak a, M))}(M/T_R(\frak a, M)).$$  Consequently, we can (and do) assume that $T_R(\frak a, M)=0,$ and with this assumption our aim is to show that $\Ann_R(H_{\frak a}^c(M))= \Ann_R(M).$ To this end, as $$\Ann_R(M)\subseteq \Ann_R(H_{\frak a}^c(M)),$$ it is enough for us to prove that $$\Ann_{R/\Ann_R(M)}(H_{\frak a}^c(M))=0.$$
Since by Lemma 2.1, ${\rm cd}(\frak a, R/\Ann_R(M))=\dim R/\Ann_R(M)= c$,  it follows that it is enough for us to show that
 $$\Ann_{R/\Ann_R(M)}( H^c_{\frak a(R/\Ann_R(M))}(M))=0.$$

 We can, and do,  assume henceforth in this proof that ${\rm cd}(\frak a, R)=c=\dim R$ and  that $M$ is a faithful $R$-module. Hence it is sufficient for us to show that $\Ann_R(H_{\frak a}^c(M))=0.$ To this end, let $x\in \Ann_R(H_{\frak a}^c(M)).$  Our strategy is to show that $H_{\frak a}^c(xM)=0.$  To do this,
it is sufficient for us to show that  $H_{\frak aR_{\frak p}}^c(xM_{\frak p})=0,$ for all  $\frak p\in \Spec R$. If  $\dim_{R_{\frak p}}(xM_{\frak p})< c,$ then
the assertion follows from Grothendieck's vanishing theorem (see \cite[Theorem 6.1.2]{BS}).  Hence we may assume that $\dim_{R_{\frak p}}(xM_{\frak p})= c.$ Then  $\dim_{R_{\frak p}}M_{\frak p}= c.$ Now, if  $H_{\frak aR_{\frak p}}^c(M_{\frak p})=0,$ then it follows from Lemma {\rm2.1}  that
$H_{\frak aR_{\frak p}}^c(xM_{\frak p})=0$, and there is nothing to prove. We therefore make the additional assumption that $H_{\frak aR_{\frak p}}^c(M_{\frak p})\neq 0.$ Now, as $xH_{\frak aR_{\frak p}}^c(M_{\frak p})= 0$, it follows from \cite[Corollary 3.11]{ASN} that
 $H_{\frak aR_{\frak p}}^c(xM_{\frak p})= 0.$  Consequently,  $H_{\frak a}^c(xM)= 0.$ Hence ${\rm cd}(\frak a, xM) <c$ and since $T_R(\frak a, M)=0,$ we deduce that $xM=0.$ Since $M$ is faithful, it follows that $x=0,$ as required.\qed\\

The first application of Theorem 2.3 extends the main result of \cite[Theorem 2.6]{BAG}.\\

\begin{cor}{\rm (cf. \cite[Theorem 2.6]{BAG})}
Let $(R, \frak m)$ be a local (Noetherian) ring and $M$ a finitely generated $R$-module.
Then $\Ann_R(H_{\frak m}^{\dim M}(M))= \Ann_R(M/T_R(M)).$
\end{cor}
\proof The result follows from Theorem 2.3 and the standard fact ${\rm cd}(\frak m, M)=\dim M$. \qed\\

\begin{cor}
Let $R$ be a Noetherian ring with finite dimension $c$ and $\frak a$ an ideal of $R$ such that ${\rm cd}(\frak a, R)=c$. Then $$\Ann_R(H_{\frak a}^c(R))= T_R(\frak a, R).$$
\end{cor}
\proof The assertion follows from Theorem 2.3\qed\\

\begin{rem}
Let $R$ be a Noetherian ring, $\frak a$ an ideal of $R$ and $M$ a non-zero finitely generated $R$-module with finite cohomological dimension $c:={\rm cd}(\frak a, M)$. Let $\{M_i \}_{i=0}^c$ be a filtration of submodules of $M$ such that  for each integer $0\leq i \leq c$, $M_i$ is the largest submodule of $M$ with ${\rm cd}(\frak a, M_i) \leq i$. Then $T_R(\frak a,M)=M_{c-1}$ and by \cite[Proposition 2.3]{ASN}, we have $$T_R(\frak a,M)=H_{\frak b}^0(M)=\cap_{{\rm cd}(\frak a, R/{\frak p_j})=c}N_j,$$ where $0=\cap_{j=1}^nN_j$ denotes a reduced primary decomposition of the zero submodule $0$ in $M$, $N_j$ is a $\frak p_j$-primary submodule of $M$ and $\frak b=\Pi_{{\rm cd}(\frak a, R/{\frak p_j})\neq c}\frak p_j.$
\end{rem}

\begin{cor}
Let $R$ be a Noetherian ring and $\frak a$ an ideal of $R$. Let $M$ be a non-zero finitely generated $R$-module with finite dimension $c$ such that ${\rm cd}(\frak a, M)=c$. Then $$\Ann_R(H_{\frak a}^c(M))=\Ann_R(M/H_{\frak b}^{0}(M))=\Ann_R(M/\cap_{{\rm cd}(\frak a, R/{\frak p_j})=c}N_j).$$ Here $0=\cap_{j=1}^nN_j$ denotes a reduced primary decomposition of the  zero submodule $0$ in $M$ and $N_j$ is a $\frak p_j$-primary submodule of $M$, for all $j=1,\dots ,n$ and $\frak b=\Pi_{{\rm cd}(\frak a, R/{\frak p_j})\neq c}\frak p_j.$
\end{cor}
\proof The assertion follows from Theorem {\rm2.3} and  Remark 2.6.\qed\\
\begin{cor}
Let $R$ be a Noetherian ring with finite dimension $c$ and $\frak a$ an ideal of $R$ such that ${\rm cd}(\frak a, R)=c$. Then $$\Ann_R(H_{\frak a}^{c}(R))=H_{\frak b}^{0}(R)=\cap_{{\rm cd}(\frak a, R/{\frak p_j})=c}\frak q_j,$$ where $0=\cap_{j=1}^n\frak q_j$ is a reduced primary decomposition of the zero ideal of $R$, $\frak q_j$ is a $\frak p_j$-primary ideal of $R$, for all $1\leq j\leq n$ and $\frak b=\Pi_{{\rm cd}(\frak a, R/{\frak p_j})\neq c}\frak p_j.$
\end{cor}
\proof The result follows from Corollary 2.7.\qed\\
 \begin{cor}
Let $R$ be a Noetherian ring with finite dimension $c$ and $\frak a$ an ideal of $R$ such that ${\rm cd}(\frak a, R)=c$. Then the following conditions are equivalent:

$\rm(i)$ $\Ann_R(H_{\frak a}^c(R))=0$.

$\rm(ii)$ $\Ass_R R = \{\frak p\in \Spec R|\,{\rm cd}(\frak a, R/\frak p)=c\}$.
\end{cor}
\proof The result follows from Corollary 2.8.\qed\\

\begin{cor}
Let $R$ be a Noetherian domain with finite dimension $c$ and $\frak a$ an ideal of $R$ such that ${\rm cd}(\frak a, R)=c$. Then $\Ann_R(H_{\frak a}^c(R))=0$.
\end{cor}
\proof Since $\Ass_RR={0}$, the assertion follows immediately from Corollary 2.9.\qed\\

\begin{cor}
Let $(R, \frak m)$ be a local (Noetherian) ring and $\frak a$ an ideal of $R$ such that ${\rm grade}\, \frak a={\rm cd}(\frak a, R)$. Then $\Ass_R R \subseteq \{\frak p\in \Spec R|\,{\rm cd}(\frak a, R/\frak p)={\rm cd}(\frak a, R)\}$.
\end{cor}
\proof Since ${\rm grade}\, \frak a={\rm cd}(\frak a, R)$, it follows from \cite[Theorem 3.3]{Lyn} that $\Ann_R(H_{\frak a}^c(R))=0$.
Moreover, in view of the proof of Theorem 2.3, $$ H_{\frak a}^c(R)\cong H_{\frak a}^c(R/T_R(\frak a, R)).$$
Hence $T_R(\frak a, R)\subseteq \Ann_R(H_{\frak a}^c(R))$, and so $T_R(\frak a, R)=0$. Now, the assertion follows from  Remark 2.6.\qed\\

\begin{cor}
Let $R$ be a Noetherian ring and $\frak a$ an ideal of $R$. Let $M$ be a finitely generated $R$-module with finite dimension $c$ such that ${\rm cd}(\frak a, M)=c$ and $x\in R$. Then $H_{\frak a}^c(xM)=0$ if and only if $xH_{\frak a}^c(M)=0$. In particular, $\Ann_RH_{\frak a}^c(M)=0$ if and only if ${\rm cd}(\frak a, rM)={\rm cd}(\frak a, M)$, for every non-zero element $r$ of $R$.
\end{cor}
\proof The assertion follows from Theorem {\rm2.3}.\qed\\
\begin{cor}
Let $R$ be a Noetherian ring and $\frak a$ an ideal of $R$. Let $M$ be a finitely generated $R$-module with finite dimension $c$ such that ${\rm cd}(\frak a, M)=c$. Then \\

${\rm(i)}$  $\Ann_R(H^c_{\frak a}(M))=\Ann_R(M$), whenever $\Ass_RM\subseteq \{\frak p\in \Supp M|\,{\rm cd}(\frak a, R/\frak p)=c\}$.\\

${\rm(ii)}$ $ \Rad(\Ann_R(H^c_{\frak a}(M)))=\cap_{\frak p\in \Ass_RM, \,{\rm cd}(\frak a, R/\frak p)=c}\frak p=\cap_{\frak p\in \Ass_R(M/T_R(\frak a, M))}\frak p.$\\

${\rm(iii)}$ $V(\Ann_R(H^c_{\frak a}(M)))=\Supp(M/T_R(\frak a, M)).$
\end{cor}
\proof
(i) follows from Corollary 2.7. To show (ii), use   \cite[Proposition 2.6(iii)]{ASN} and Corollary 2.7.
In order to prove (iii), in view of Theorem 2.3 we have $$V(\Ann_RH^c_{\frak a}(M))=V(\Ann_RM/T_R(\frak a, M))=\Supp (M/T_R(\frak a, M)).$$\qed

\section{Attached primes of local cohomology modules}

It will be shown in this section that the subjects of the previous
section can be used to investigate the attached prime ideals of local
cohomology modules. In fact, we will generalize and improve the main result of
Hellus (cf. \cite[Theorem 2.3]{He}). The main result
is Theorem 3.7. The following proposition will serve to shorten the proof of the main theorem. We begin with\\

\begin{defn}
 Let $L$ be an $R$-module. We say that a prime ideal  $\frak p$ of $R$ is an {\it  attached prime of}  $L$, if there exists a submodule $K$
 of $L$ such that $\frak p=\Ann_R(L/K)$. We denote by $\Att_R L$ the set of attached primes of $L$.
\end{defn}

When $M$ is {\it representable} in the sense of \cite{Ma} (e.g. Artinian or injective), our definition of $\Att_R L$ coincides with that of Macdonald, Sharp \cite{Ma}, \cite{Sh}.

It follows easily from the definition that, if $\frak p\in \Att_R L$, then $L\otimes_R R/\frak p\neq0$.  This is used in the proof of Theorem 3.3.

\begin{lem}
Let $R$ be a Noetherian ring and $L$ an $R$-module. Then, the set of minimal elements of $V(\Ann_R(L))$ coincides with that of
$\Att_RL$. In particular, $$\Rad(\Ann_R(L))= \cap_{\frak p\in \Att_R L}\frak p.$$
\end{lem}

\proof The assertion follows from the Definition 3.1 and the fact that $\frak p=\Ann_R(L/\frak p L)$ for each minimal prime $\frak p$ over $\Ann_R(L)$. \qed\\

We are now ready to state and prove the first main result of this section, that gives us an upper bound for the attached primes of $\Att_R H_{\frak a}^{{\rm cd}(\frak a, M)}(M).$ Before this, we note that as $\Ann_R(M)\subseteq\Ann_R(H^i_{\frak a}(M)),$ it follows that $\Att_R(H^i_{\frak a}(M))\subseteq \Supp M,$ for every finitely generated $R$-module $M.$\\

\begin{thm}
Let $R$ be a Noetherian ring and $\frak a$ an ideal of $R$. Let $M$ be a non-zero finitely generated $R$-module such that $c:={\rm cd}(\frak a, M)$ is finite. Then
$$\Att_RH^c_{\frak a}(M)\subseteq \{\frak p\in \Supp M|\,{\rm cd}(\frak a, R/\frak p)= c\}.$$
\end{thm}
\proof As $\Att_RH_{\frak a}^c(M)\subseteq \Supp M$, it follows from Lemma {\rm2.1}  that
 $$\Att_RH^c_{\frak a}(M)\subseteq \{\frak p\in \Supp M|\,{\rm cd}(\frak a, R/\frak p)\leq c\}.$$
  Now, it is enough to show that
$$\Att_RH^c_{\frak a}(M)\subseteq \{\frak p\in \Supp M|\,{\rm cd}(\frak a, R/\frak p)\geq c\}.$$
 To this end, let $\frak p\in \Att_RH^c_{\frak a}(M)$. Then   $$\frak p/\Ann_R(M)\in \Att_{R/\Ann_R(M)}H^c_{\frak a}(M),$$ and so
$$H^c_{\frak a}(M)\otimes_{R/\Ann_R(M)}R/\frak p\neq 0.$$ Now, as $$H^c_{\frak a}(M)\cong H^c_{\frak a}(R/\Ann_R(M))\otimes_{R/\Ann_R(M)}M,$$ it follows that $$H^c_{\frak a}(R/\Ann_R(M))\otimes_{R/\Ann_R(M)}M\otimes_{R/\Ann_R(M)}R/\frak p\neq 0.$$
Consequently $$H^c_{\frak a}(R/\frak p)\otimes_{R/\Ann_R(M)}M\neq 0,$$  and thus $ H^c_{\frak a}(R/\frak p)\neq 0$, as required.\qed\\

 The next corollary reproves the main result of \cite[Theorem 2.5]{Di}.\\

\begin{cor}
 Let $R$ be a Noetherian ring and $\frak a$ an ideal of $R$. Let $M$ be a non-zero finitely generated $R$-module of finite dimension $d$.
Then
 $$\Att_RH^d_{\frak a}(M)=\{\frak p\in \mAss _R M|\,{\rm cd}(\frak a, R/\frak p)= d\}.$$
\end{cor}
\proof We can (and do) assume that $H^{d}_{\frak a}(M)\neq 0$. Then in view of Theorem {\rm3.3}, we have $$\Att_RH^d_{\frak a}(M)\subseteq\{\frak p\in \mAss _R M|\,{\rm cd}(\frak a, R/\frak p)= d\}.$$
Now, in order to show the inverse containment, let $\frak p$ be a prime ideal of $R$ such that  $\frak p\in \mAss _R M$ and  ${\rm cd}(\frak a, R/\frak p)= d$.
Then in view of \cite[Proposition 2.6]{ASN}, we have $\frak p\in \Ass_R M/T_R(\frak a, M)$, and so by Corollary {\rm2.13(iii)},
$\frak p\in V(\Ann_R(H^d_{\frak a}(M)))$. Hence in view of  Lemma {\rm3.2},  $\frak p\in \Att_RH^d_{\frak a}(M)$, as required. \qed\\

\begin{cor}
Let $R$ be a Noetherian ring and $\frak a$ an ideal of $R$. Let $M$ be a finitely generated $R$-module such that $c:={\rm cd}(\frak a, R)$ is finite. Then
\begin{eqnarray*}
\{\frak p\in \Spec R|\, {\rm cd}(\frak a, R/\frak p)=c=\dim R/\frak p\} \subseteq \Att_RH^c_{\frak a}(R) \subseteq \{\frak p\in \Spec R|\,{\rm cd}(\frak a, R/\frak p)=c\}.
\end{eqnarray*}
Moreover,  if $c=\dim R$, then  $\Supp H^c_{\frak a}(R)\subseteq V(\frak a+T_R(\frak a, R)).$
\end{cor}

\proof
In order to prove the first containment, in view of Theorem 3.3, it is enough to show that
$$\{\frak p\in \Spec R|\,{\rm cd}(\frak a, R/\frak p)=c=\dim R/\frak p\}\subseteq \Att_RH^c_{\frak a}(R).$$
To this end, let $\frak p\in \Spec R$ be such that ${\rm cd}(\frak a, R/\frak p)=c=\dim R/\frak p.$ Then in view of Corollary 3.4, we have $\frak p\in \Att_RH^c_{\frak a}(R/\frak p)$. Now, from the exact sequence
$$0 \longrightarrow \frak p \longrightarrow R \longrightarrow R/{\frak p} \longrightarrow 0,$$
and the right exactness of $H^c_{\frak a}(\cdot)$, we deduce that $\frak p \in \Att_RH^c_{\frak a}(R)$, as required.

In addition, in order to show the last inclusion,  use Corollary 2.5 and the fact that $\Supp H^c_{\frak a}(R)\subseteq V(\frak a).$ \qed\\

\begin{lem}
Let $R$ be a Noetherian domain of finite dimension $d$ and $\frak a$ an ideal of $R$ such that ${\rm cd}(\frak a, R)=d-1$. Then
$\Ann_R(H^{d-1}_{\frak a}(R))=0.$
\end{lem}

\proof Set $J:=\Ann_R(H^{d-1}_{\frak a}(R))$. We now suppose that $J\neq 0$,  and look for a contadiction. To this end, from
 the right exactness of the functor $H^{d-1}_{\frak a}(\cdot)$, we deduce that
 $H^{d-1}_{\frak a}(R)\cong H^{d-1}_{\frak a}(R/J^t)$,  for all integers $t$,  and  so by Lemma {\rm2.1} we have ${\rm cd}(\frak a, R/J^t)=d-1$.
 Now, as $R$ is integral domain and $J^t\neq 0$, it follows that $\dim R/J^t=d-1$, for all integers $t$. Moreover, in view of
 Corollary {\rm3.4},  there is a $\frak p\in \mAss_RR/J^t$  such that ${\rm cd}(\frak a, R/\frak p)=d-1$, for all integers $t$.
  Now, let $\frak q_t$  be the $\frak p$-primary component of $J^t$.  Then, in view of Corollary {\rm2.7},  we have
 $$J=\Ann_R(H^{d-1}_{\frak a}(R/J^t))\subseteq \frak q _t, $$  for all integers $t$. Consequently, we obtain that
 \begin{center}
 $JR_{\frak p}\subseteq \bigcap_{t\geq0}\frak q_tR_{\frak p}=\bigcap_{t\geq0} J^tR_{\frak p},$
 \end{center}
 and so the  Krull's Intersection Theorem implies that  $JR_{\frak p}=0$. َAs $R$ is an integral domain, it follows that $J=0$, which is a contradiction.   \qed\\

Now we are prepared to prove the second main theorem of this section, which
is a generalization of the main result of Hellus (cf. \cite[Theorem 2.3]{He}).\\

\begin{thm}
Let $R$ be a Noetherian ring of finite dimension $d$ and $\frak a$ an ideal of $R$ such that $H^d_{\frak a}(R)=0$. Then

$$\Att_RH^{d-1}_{\frak a}(R)= \{\frak p\in \Spec R|\,{\rm cd}(\frak a, R/\frak p)=d-1\}.$$
\end{thm}
\proof We can (and do) assume that $H^{d-1}_{\frak a}(R)\neq 0$. Then in view of Theorem {\rm3.3}, we have $$\Att_RH^{d-1}_{\frak a}(R)\subseteq \{\frak p\in \Spec R|\,{\rm cd}(\frak a, R/\frak p)=d-1\}.$$ Now, let $\frak p$ be a prime ideal of $R$ such that ${\rm cd}(\frak a, R/\frak p)=d-1$,
and so $\dim R/\frak p\geq d-1$.
Hence, using  Corollary {\rm2.10} and Lemma {\rm3.6} we can easily see that $\Ann_{R/\frak p}H^{d-1}_{\frak a}(R/\frak p)=0.$ Accordingly, we have $\Ann_RH^{d-1}_{\frak a}(R/\frak p)=\frak p.$ Thus,  in view of the definition,  we have $\frak p\in\Att_RH^{d-1}_{\frak a}(R/\frak p).$ Therefore, from the exact sequence
 $$0 \longrightarrow \frak p \longrightarrow R \longrightarrow R/{\frak p} \longrightarrow 0,$$
 and the right exactness of the functor $H^{d-1}_{\frak a}(\cdot)$, we deduce that $\frak p \in \Att_RH^{d-1}_{\frak a}(R)$, as required.\qed\\

The following lemma, which is a consequence of the Lichtenbaum-Hartshorne vanishing theorem, is assistant in the proof of Theorem {\rm3.9}.\\

\begin{lem}
Let $(R, \frak m)$ be a complete local (Noetherian) ring and $\frak a$ an ideal of $R.$ Let $\frak p$ be a prime ideal of $R.$ Then
${\rm cd}(\frak a, R/\frak p)=\dim R/\frak p$ if and only if $\Rad(\frak a+\frak p)=\frak m.$
\end{lem}
\proof  Let ${\rm cd}(\frak a, R/\frak p)=\dim R/\frak p:=c.$ Then $H^c_{\frak a}(R/\frak p)\neq0$, and so according to the Lichtenbaum-Hartshorne vanishing theorem (see \cite[Theorem 8.2.1]{BS}), $\Rad(\frak a+\frak p)=\frak m$. \\
In order to prove of the opposite direction, use \cite[Theorems 4.2.1 and 6.1.4]{BS}. \qed\\

The next theorem which is a consequence of Theorem 3.7, improves the main result of \cite[Theorem 2.3]{He}. If  $(R, \frak m)$ is a local  ring,
then we use $\hat{R}$ to denote the completion of $R$ with respect to the $\frak m$-adic topology. \\

\begin{thm}
Let $(R, \frak m)$ be a local (Noetherian) ring of dimension $d$. Let $\frak a$ be an ideal of $R$ such that $\dim R/\frak a=1$ and
$H^d_{\frak a}(R)=0$.  Then
$${\rm Assh}_RR\subseteq \Att_RH^{d-1}_{\frak a}(R).$$
If, in addition, $R$ is complete, then
$$\Att_RH^{d-1}_{\frak a}(R)=\{\frak p\in \Spec R|\, \dim R/\frak p=d-1\,\text{and}\, \Rad(\frak a+\frak p)=\frak m\}\cup {\rm Assh}_RR.$$
\end{thm}
\proof  For the proof of ${\rm Assh}_RR\subseteq \Att_RH^{d-1}_{\frak a}(R)$, let $\frak p\in{\rm Assh}_RR$. Then, in view of Theorem {\rm3.7}, it is enough
to show that ${\rm cd}(\frak a, R/\frak p)=d-1$.  To do this, as 
${\rm Assh}_RR=\{\frak q\cap R|\, \frak q\in {\rm Assh}_{\hat{R}}\hat{R}\}$, $\dim \hat{R}/\frak a\hat{R}=1$ and
 ${\rm cd}(\frak a\hat{R}, \hat{R}/\frak p\hat{R})={\rm cd}(\frak a, R/\frak p)$, by using Lemma 2.1  without loss of generality, we may assume that $R$ is complete.
 Now, since $H^d_{\frak a}(R)=0$, it follows that $H^d_{\frak a}(R/\frak p)=0$. Thus $\frak a+\frak p$ is not $\frak m$-primary, and so
 $\dim R/(\frak a+\frak p)=1$ (note that $\dim R/\frak a=1$). Let $\frak q$ be a minimal prime over $\frak a+\frak p$ such that $\dim R/\frak q=1$.
 Then, it is easy to see that $\frak q$ is also  minimal over $\frak a$. Next, as $R$ is catenary, it yields that
 $\dim R_{\frak q}/\frak p R_{\frak q}=d-1$, and so we have
 $$(H_{\frak a}^{d-1}(R/\frak p))_{\frak q}\cong H_{\frak q R_{\frak q}}^{d-1}(R_{\frak q}/\frak pR_{\frak q})\neq0.$$
 Therefore $H_{\frak a}^{d-1}(R/\frak p)\neq0$, and hence ${\rm cd}(\frak a, R/\frak p)=d-1$; so that
 ${\rm Assh}_RR\subseteq \Att_RH^{d-1}_{\frak a}(R)$. Consequently, Theorem {\rm3.7} enables us to deduce that $\Att_RH^{d-1}_{\frak a}(R)$ equals
 the set
 $$\{\frak p\in \Spec R|\, {\rm cd}(\frak a, R/\frak p)=\dim R/\frak p=d-1\}\cup\{\frak p\in {\rm Assh}_RR|\,{\rm cd}(\frak a, R/\frak p)=d-1\}.$$
 Now, using Lemma {\rm3.8}, we see that $\Att_RH^{d-1}_{\frak a}(R)$ is equal with the set
  $$\{\frak p\in \Spec R|\, \dim R/\frak p=d-1\,\text{and}\, \Rad(\frak a+\frak p)=\frak m\}\cup \{\frak p\in {\rm Assh}_RR|\,{\rm cd}(\frak a, R/\frak p)=d-1\}.$$
  Therefore, it follows from  Theorem {\rm3.7} and ${\rm Assh}_RR\subseteq \Att_RH^{d-1}_{\frak a}(R)$ that
  $$\Att_RH^{d-1}_{\frak a}(R)=\{\frak p\in \Spec R|\, \dim R/\frak p=d-1\,\text{and}\, \Rad(\frak a+\frak p)=\frak m\}\cup {\rm Assh}_RR,$$
  as required.  \qed\\

\begin{rem}
As a main result, it has been proved in \cite[Theorem 2.3]{He}, if $(R, \frak m)$ is a complete local ring of dimension $d$ and $\frak a$ an ideal of $R$ such that $\dim R/\frak a=1$ and $H^d_{\frak a}(R)=0$, then
$$\Att_RH^{d-1}_{\frak a}(R)=\{\frak p\in \Spec R|\, \dim R/\frak p=d-1\,\text{and}\, \Rad(\frak a+\frak p)=\frak m\}\cup {\rm Assh}_RR.$$
The proof of \cite[Theorem 2.3]{He} relies heavily on \cite[Theorems 2.4 and 2.5]{HS}, and these results are not true in the case that $d=1$.
Indeed, if $(R, \frak m)$ is a complete local domain of dimension $d=1$, then it is easy to see that these results are not true. Thus
\cite[Theorem 2.3]{He} needs correction, nevertheless their proofs are valid in the case that $d\geq2$. However, Theorem {\rm3.9} recovers
the corrected version of \cite[Theorem 2.3]{He}.
\end{rem}

We end the paper with the following question:\\

{\bf Question.} $\rm(i)$  We have shown in Theorem {\rm3.3} that for an ideal $\frak a$  of a Noetherian ring $R$ and a finitely generated $R$-module $M$ of finite cohomological dimension $c$,
$$\Att_RH^c_{\frak a}(M)\subseteq\{\frak p\in \Supp M|\,{\rm cd}(\frak a, R/\frak p)= c\}.$$
Is the above containment equality?

$\rm(ii)$ In Theorem 2.3 we have determined $\Ann_R(H_{\frak a}^{\dim M}(M))$. Is it possible to determine 
$\Ann_R(H_{\frak a}^{{\rm cd}(\frak a, M)}(M)) ?$

\begin{center}
{\bf Acknowledgments}
\end{center}
The authors are deeply grateful to the
referee for his or her valuable suggestions on the paper and for
drawing the authors'  attention to Theorem 2.3 and Remark  2.6.  Also, we would like to thank Prof. Kamal Bahmanpour for useful discussions, and also the Institute for Research in Fundamental Sciences (IPM) for the financial support.

\end{document}